\documentclass{amsart}
\newtheorem{thm}{Theorem}[section]

\newtheorem{lemma}[thm]{Lemma}
\newtheorem{prop}[thm]{Proposition}
\theoremstyle{definition}

\theoremstyle{remark}

\numberwithin{equation}{section}

\usepackage{amssymb}
\usepackage{graphicx}
\usepackage{epstopdf}

\newcommand{\dd}[2]
{
{{\partial #1}   \over {\partial #2}}
}

\newcommand{\ddd}[2]
{
{{d #1}   \over {d #2}}
}

\newcommand{\ddt}{{d \over dt}}

\begin{document}

{\bf \large The only syzygy-free solution is Lagrange's.}
 
 \smallskip  
{\small by  Richard Montgomery,   Mathematics Dept., UC Santa Cruz,
   Santa Cruz, CA, 95064, USA,  supported in part by NSF grant DMS-20030177.}
\medskip 

{ \small  {\bf Absract.}
A syzygy in the three-body problem is a collinear instant. 
We prove that with the exception of
Lagrange's solution every solution to  the zero angular momentum Newtonian three-body problem suffers syzygies.  The proof works    for all  mass
ratios.  
}

\vskip .3cm

We consider the Newtonian three-body problem
with   zero angular momentum  and negative energy. Masses 
are positive, but     arbitrary. 
A `syzygy'   means  an eclipse:  an instant at which 
the three masses are collinear.  
\begin{thm}. 
\label{thm:main}  Every  solution  admits a syzygy 
except one:  the   Lagrange homothety solution.
\end{thm}
 Solutions are defined   
over their    maximal interval of existence
and   analytically  continued through   binary collisions a la 
  Levi-Civita  \cite{Levi}. Binary collisions   counts as syzygies. 
A solution   cannot be  
   extended past a finite time $t = b$ if and only if 
    as $t \to b$ the three  positions of the  three bodies 
    tend to the same point.    In other words, a solution   fails  to exist past a certain time
  if and only  if it ends in triple collision at that time.  (See \cite{Sundman},
    \cite{Moeckel} or \cite{Saari} )

The Lagrange homothety  solution \cite{Lagrange} begins and ends in triple collision.
At every  other instant of its existence the masses form an equilateral triangle.  This triangle evolves
by  homothety (scaling). Half-way through
its evolution the three   bodies are instantaneously at rest,  forming
an equilateral   triangle whose size is  determined by the value of the negative energy.   
 
 In \cite{Mont}  I proved   theorem \ref{thm:main} upon imposing two  additional 
  hypotheses on solutions: that they are  bounded,
 and that  they do  
 not end in triple collision.       The contribution
 of the present paper is to dispense with these   hypotheses.

 I first dispense with the hypothesis on collision, keeping the boundedness
 hypothesis.  Again, in  \cite{Mont}  I  proved that
  bounded solutions which do not end in collision have syzygies. 
  The same proof, plus   invariance  of the
 equations and zero angular momentum condition under time reversal
proves existence of  syzygies for solutions which are bounded and do not
 begin in triple collision.   All that remains of the
 bounded soltuions are those,  excluding Lagrange,  which begin and end in triple collision.     The proof for these solutions will 
 follow  the same   qualitative lines as  \cite{Mont}.  
  According to    Moeckel [1989]  (see the corollary at the top of p. 53)    there are,
    for generic mass ratios,      an infinite number of these finite-interval solutions bi-asymptotic   to triple collision. 
  
   Moeckel,  Chenciner and others have  pointed out  that dispensing with
   the boundedness hypothesis on solutions  ought to be easy.  In unbounded negative energy solutions  two of the masses must  form a bound pair
   with the third mass far away 
   for long periods of time. During these long times the bound 
    pair moves according to  a differential equation which
   is a slight, but time-dependent perturbation of the Kepler equation and so the
   pair should spin about each other frequently crossing the line joining their
   center of mass to the distant mass, and thus making syzygies.  However
   I was unable to turn this idea  into a proof.  The difficulties include 
   the existence of  oscillatory unbounded solutions, and  the 
   difficulty of establishing syzygies for systems looking like 
   highly eccentric nearly  Keplerian orbits subject to small time-dependent perturbations
   concentrated along the semi-major axis of the orbits.  
   Instead, I  use the  methods of \cite{Mont}.  The    bulk of this paper
     is devoted to proving the existence of infinitely many syzygies for unbounded solutions with zero angular momentum.  I expect  a
     more skilled analyst could get a more  direct proof based on the Kepler idea,
     and valid for  unbound negative energy solutions
     with nonzero angular momentum.

  {\bf   Motivation.}    I have
    been trying for some time to establish a   symbolical dynamical description for the zero angular momentum three-body.   The symbols are to be the  syzygies, marked as 1, 2 , 3, depending on which mass crosses between the other two.   See \cite{pants},
    \cite{braids}.   I have 
    successfully established   a 
    complete symbolical dynamical description  if I allow myself to  change the potential  
    from the Newtonian  $1/r$ potential to  the $1/r^2$ potential 
    and if I take   the three  masses to be equal.     
  Theorem  \ref{thm:main}     is  a first step toward  the more interesting Newtonian case.
    The theorem  asserts that, with one exception , every solution has  syzygies, 
    and hence a  symbol sequence.    
    
        \medskip 
  {\bf Proof. } 
    
    We continue to use the methods of   \cite{Mont} where we introduced
    the  ``height'' variable $z$  on  the three-body   configuration space minus
    triple collisions.
    The   crucial properties of this $z$ are
     \begin{equation}
     -1 \le z \le 1 
      \end{equation}
    \begin{equation}|z| =1 \iff \hbox{equilateral } 
     \end{equation}
    \begin{equation} z = 0    \iff \hbox { syzygy }  
     \end{equation}
    and that along any solution 
  \begin{equation}d(f \dot z )/dt   =  -q z,  \hskip 1cm f > 0, g \ge 0
  \label{eqn:DE}
  \end{equation}
  where $f$ is a  smooth function on shape space, $q$ is a
  smooth  function on the tangent space
  to   shape space, 
  \begin{equation} 
  \label{eqn:Lagcase} 
  q = 0 \iff  \hbox {tangent to Lagrange homothety} .
   \end{equation}
  and
  \begin{equation} f \to \infty  \iff  \hbox {unbounded}.
  \label{eqn:fbd}
  \end{equation}
We recall that  a solution is called bounded if
  all the distances $r_{ij}$ between the pairs $i,j$ of masses
  are bounded functions of time.  Thus the solution is   unbounded if
  $\lim \sup r_{ij} (t) = + \infty$ for some mass pair $i,j$. 

{\bf The bounded case.}  
Let $x$ be a solution as per the theorem, and suppose
it to be bounded.  Thus $x$ is a  bounded  zero-angular momentum
negative energy  solution to the three-body problem  besides the Lagrange solution.
We may suppose it is not  a  collinear solution since
every instant of a collinear solution is a syzygy. 
Reflecting a solution about a line effects the transformation $z \to -z$,
and a time reflection   $t \to -t$ effects  the transformation $\dot z \to - \dot z$.
Using these symmetries and time translation , we may assume at some initial time $t= -\epsilon$ we have $z > 0$ and $\dot z \le 0$.  
  Because the solution  is not the Lagrange homothety solution, 
  $z$ cannot be identically $1$ and  
  $q$ must be positive along the solution (see  (\ref{eqn:Lagcase})),    
  It follows from (\ref{eqn:DE})  that    $(d/dt)(f \dot z) < 0$.  
  In particular $\dot z =0$ identically is impossible.  Upon translating  time   
    forward slightly from $-\epsilon$ to $0$ 
  we will have  
    $\dot z < 0$.  Now  we have  $z(0) > 0$
    and $\dot z (0) < 0$.  We must prove that at   some   finite time $t =b$ later  
    we have   $z(b) = 0$.

 According to (\ref{eqn:DE})  $f \dot z$ is strictly decreasing as long as $z > 0$.
 Since    $f$ is positive,
  the derivative   $\dot z$ must remain negative 
  over any interval $[0,b)$ of time during which $z(t) > 0$.
  Thus  $z (t) $ is monotonically decreasing over every interval of time  
   $[0, b)$ for which the solution exists and for which   $z(t) >0$.
      The solution cannot
  fail to exist in   such interval, because the only way it can terminate itself 
  is by ending in triple collision.  But non-collinear solutions
  which end in triple collision must asymptote to the Lagrange solution 
 \cite{Moeckel} \footnote{See Table 1 there, the entry $dim(St (R)$ with 
 $R = C^*$.  The linearization at $C^*$ for collinear motion also has   $dim(St (R) = 1$,
 showing that the stable manifold ingoing to a collinear triple collision $C^*$ lies within
 the collinear submanifold. }  implying  $z \to 1$ or $z \to -1$ which
  we have excluded.   
  Hence either $b < \infty$ and  $z(b) = 0$, in which case  we have our
   syzygy, or $b = \infty$
  and the solution stays in the upper hemisphere $z > 0$ for all positive time.
 We   invoke the hypothesis that the solution is bounded
 to   exclude the second possibiliy.
 
 So, suppose that $z > 0$ on $[0, \infty)$, that $\dot z (0) < 0$
 and that the motion is bounded.  According to the bound (\ref{eqn:fbd})
  the function $1/f$    is bounded away from zero along our solution, so that    $1/f  > k$ on $[0, \infty)$ for some positive constant  $k$.  
  Now $\dot z (0)$ is negative by assumption, and $f(0)$ is positive
  so that $f(0) \dot z (0)  = - a < 0$ is negative.
      According to the  differential equation (\ref{eqn:DE})
  $f \dot z $ is monotonically decreasing so that $f (t) \dot z(t) < -a$.
  Then $\dot z = {1 \over f} (f \dot z) < -ka$.
  But then  
  $$z(t) = z(0) + \int_{0} ^t  \dot  z  dt < z(0) - k a t .$$
  which   violates the positivity of $z$ as soon as $t > z(0)/ka$.
  This contradiction shows that in fact $z$ has a zero
  before time $t = z(0)/ka$.  
  
  {\bf }    
    
  {\bf The unbounded case.}

 There are two types of unbounded solutions, 
  escape and   oscillatory.
A solution is an escape solution if 
$\lim r_{ij} (t) = + \infty$ for some pair $ij$.
It is an oscillatory solution if for some pair 
$\lim \sup r_{ij} = \infty$ while for every pair
$\lim \inf r_{ij} < \infty$.  The existence of oscillatory type
unbounded solutions was established by Sitnikov \cite{Sitnikov}.
Our proof   deals simultaneously with both types. 
%yes?? 

The function $\max_{ij}  r_{ij}$ is a measure of the
size of the configuration.  Another equivalent measure
is $R$ where   $R^2 = I$
is the moment of inertia: 
 $I = \Sigma m_i m_j r_{ij} ^2 / \Sigma m_i$.  The $m_i$ are the values of the masses. 
Then 
$$ c \max_{ij}  r_{ij} < R < C \max_{ij}  r_{ij} .$$
{\bf where here and throughout $c, C$ denote positive constants depending
only on the masses and occassionally on the energy.}  The precise values of these
constants will not be important. 
It follows  that a  motion is  an escape motion if 
$\lim_{t \to \infty}  R (t) = +\infty$ and it is oscillatory if  
$\lim \sup _t R(t) = + \infty$ while $\lim \inf _t R (t) < \infty$.
The function $f$ of the basic equation (\ref{eqn:DE})
is related to $R^2$ by
\begin{equation}
\label{eqn:f2}
f = R^2 \lambda ^2
\end{equation}
where $c \le \lambda \le C$.
Relation  (\ref{eqn:fbd}) follows   from this expression for $f$. 

Let $x$ be an  unbounded solution.   
Being unbounded, for any  $R_0 > 0$
there is a time $t$ such that $R(t) \ge R_0$.  
Rewrite the $z$ equation (\ref{eqn:DE}) by  introducing a new time variable
$\sigma$:  $$f \ddd{z}{t} = \ddd{z}{\sigma}, \hbox{ so that } 
\ddd{t}{\sigma} = f$$
The differential equation for   $z$   becomes:
${1 \over {f}} {{d^2 z} \over {d \sigma ^2}} = -q z$
which is a harmonic oscillator:  
\begin{equation}
\label{eqn:ODEosc}
 {{d^2 z} \over {d \sigma^2}} = - \omega^2 z; \omega^2 = f q 
 \end{equation}    of variable frequency $\omega$.   
If $\omega$ were to be a constant  $\omega_0$ then this would be the
equation of a    linear oscillator and the  zeros of any solution would be spaced equally
at ($\sigma$-) intervals of  length $\pi/\omega_0$.   
Returning to our case, from standard
Sturm-Liousville theory it follows that if
$\omega^2 > \omega_0 ^2$ then within   each of these  intervals of  length $\pi/{\omega_0}$ the function $z(\sigma)$
has a zero.  Let $\ell$ be the length of an interval of $\sigma$
during which $R \ge R_0$ and suppose that $\omega  \ge \omega_0$
during this interval.   In the $\sigma$ variable    escape to infinity takes  finite time
so the lengths $\ell$ will be finite.
If we can show that for $R_0$ sufficiently large
$\ell > \pi/ \omega_0$ then
we will know there is an oscillation during this interval. 
We will establish below the asymptotics:  
\begin{equation}  \label{eqn:periods}
\ell  \omega_0  \ge C R_0,  R_0 \to \infty 
 \end{equation}  
 It follows that 
 there are   many syzygies
during the   interval $\ell$  (at least  $CR_0 / \pi$ syzygies ).

The following two estimates yield  (\ref{eqn:periods})
\begin{equation}  \label{eqn:omega} \omega \ge C R_0 ^2 \end{equation}
\begin{equation} \label{eqn:interval} \ell \ge C/R_0  \end{equation}

{\bf Proving estimate \ref{eqn:omega}, the $\omega$ bound.}

Let 
$$ r = \min_{i \ne j}  r_{ij} $$
be  the minimum   distance.
 Fix the total energy $H$  to be negative and write $h = -H > 0$.
Then, as is well-known,  there exists a constant $c$ depending only on the masses
such that   the minimum intermass distance
$r$ satisfies 
\begin{equation}
\label{eqn:rbd1}
r \le c/(|H|).
\end{equation}  See equation (\ref{eqn:rbd}) of the
appendix for a proof.

The total energy is given 
$$H = (K/2) - U$$
where $K \ge 0$ is the potential energy and  
$$U = \Sigma m_i m_j /r_{ij}.$$
is the negative of the potential energy. 
Because our solution has negative energy and 
$R$ is large along intervals of the solution of it, 
we know that along any one such `long' interval that
one of the distances,  
say $r_{12}$  is much smaller
than the other two  and these other two  are  of   order $R$: 
\begin{equation}
 \label{eqn:bdsr} r_{12} = r  \le C,   r_{13}, r_{23}  \ge CR.
\end{equation}
(See the appendix , equations (\ref{eqn:rhoRbds}, \ref{eqn:distances}) for proofs.) 
Introduce  the spherical coordinates
$(R, \theta, \phi)$   used in   \cite{Mont}.   
and  the squared distance variables  
$$s_k = r_{ij} ^2$$
for $i,j,k$ any permutation of $123$.  
These systems of 
coordinates are related by    $s_k =    = R^2  \lambda (1 - \gamma_k (\theta) \cos(\phi))$
where $\gamma_k = \cos(\theta - \theta_k)$,  and $\theta = \theta_k$, 
$\phi = 0$ desribes the location of the binary collision ray  $r_{ij} = 0$.
See \cite{Mont} eq. (4.3.14).  
(The angles    $\theta_k - \theta_j$ between the collision   rays depend on the masses.)

 The function  $q$ satisfies:  
 $$q = \hbox{positive} + {-4 \cos(\phi) \over \sin(\phi)} {{\partial U} \over {\partial \phi}}.$$
 It follows that the bound (\ref{eqn:omega})  will follow from
 the bound 
\begin{equation}
 \label{eqn:bdpotl}
  {-4 \cos(\phi) \over \sin(\phi)} {{\partial U} \over {\partial \phi}} \ge C R^2 
  \end{equation}
valid for all $R$  large enough, together with the 
defining relations (\ref{eqn:ODEosc}) and (\ref{eqn:f2}). 

We proceed to establish the bound (\ref{eqn:bdpotl}). 
%For shorthand, we now set   $$c = \cos(\phi), s = \sin(\phi).$$
 We have
   \[
\begin{array}{lcl}
 {{\partial U} \over {\partial \phi}} &  =  &  
  - \Sigma {{m_i m_j} \over {r_{ij} ^3}}\dd{s_k}{\phi}   
  \end{array}
\]
and
\[
\begin{array}{lcl}  
\dd{s_k}{\phi}     &  =   & - R^2 \lambda \gamma_k  \sin(\phi)  + R^2 \dd{\lambda}{\phi} (1 - \gamma_k  (\theta) \cos(\phi))  \\
  &  =  &  \sin(\phi)  R^2 \lambda \gamma_k  + \Lambda s_k
\end{array}
\]
where $\Lambda = \dd{log \lambda}{\phi}$.  
 The function 
$\lambda$ and hence $\log \lambda$ are even functions on the shape sphere, 
where `even' and `odd' refer to behaviour under the  reflection$(\phi, \theta) \mapsto (-\phi, \theta)$  about the equator.  It follows 
 that $\Lambda$ is an odd function and so vanishes
on the equator.  Thus 
$$\Lambda = \sin(\phi) W (\phi, \theta)$$
where $W$ is a smooth function on the sphere. In particular
$W$ is uniformly bounded.    Now we have
\[
\begin{array}{lcl}  
{-1 \over {\sin(\phi)}}{1 \over {r_{ij}^3}} \dd{s_k}{\phi}     &  =    &   {1 \over {r_{ij}^3}} [ R^2 \lambda \gamma_k (\theta)  - W s_k ]  
\end{array}
\]
Since   $\gamma_3 =1$ at collision and  since  the point  in the shape sphere
representing our triangle is arbitrarily close to this  same collision point 
for $R$ big (because $r/R << 1$), we have that   $\gamma_3$ is as 
  close to $1$ as we like along our solution interval, by taking   $R$ large along the interval.  
Using this fact, and the bound   (\ref{eqn:rbd1}) we have
\[
\begin{array}{lcl}  
{-1 \over {\sin(\phi)}}{1 \over {r_{12}^3}} \dd{s_3}{\phi}   &  =   &  {1 \over {r_{12}^3}} [ R^2 C - C]  \\
& \ge &  CR^2 
\end{array}
\]
And for the other two distances   :
\[
|{-1 \over {\sin(\phi)}}{1 \over {r_{13}^3}} \dd{s_2}{\phi} |, |{-1 \over {\sin(\phi)}}{1 \over {r_{23}^3}} \dd{s_1}{\phi} |   \le C/R
 \]
Thus: 
\[
\begin{array}{lcl}
- 4 {\cos(\phi)  \over \sin(\phi)}  {{\partial U} \over {\partial \phi}} &  \ge  &   CR^2
\end{array}
\]
as claimed.

{\bf Proving the estimate on $\ell$, the bound (\ref{eqn:interval}).}
 
We will be using the   length $\rho = \| \xi |$ of the long Jacobi vector $\xi$
as a measure of escape. 
This  vector connects  the $12$ center of mass to the
   distant mass $3$.    We have
      \begin{equation}
         \label{eqn:Rbd}
         R^2 = a \rho^2 + b r_{12}^2
         \end{equation}
     where $a, b$ depend only on the masses.
 (See (\ref{eqn:R2}) of the appendix).  
       It  follows from    equation (\ref{eqn:Rbd})
        and the bound $r \le c/|H|$ on   $r = r_{12}$
     that $R$ and $\rho$ are related by
  \begin{equation}
       \label{eqn:Rrhobd}  
       c_a \rho   \le R  \le C_a \rho ,  
        \end{equation} 
         \begin{equation}
       \label{eqn:Rrhobd}  
       {1 \over C_a} R   \le \rho  \le {1 \over c_a}R,  
        \end{equation} 
        where the constants 
         $c_a, C_a$  depend only on the masses.
(These constants can  be taken   arbitrarily close  to the  constant $1/ \sqrt{a}$  by  taking $R_0$
  sufficiently large and $R > R_0$.)
  
  The desired length bound follows from the following 
  assertion
  \begin{prop}\label{prop:longsolutions}
 Let   $\rho(t)$ be  
   the length at time $t$  of the long  Jacobi vector   for a future-unbounded   negative energy solution.    Then
   %  for all $\epsilon > 0$ 
   there exists a constant $c_3$  
  such that for all $\rho_0$ sufficiently large 
   there exists a $\rho_* \ge \rho_0$ 
   %with $(1 - \epsilon) \rho_*  \ge \rho_0$
  and  two  times $t_1 < t_*$  such that 
  $\rho(t_1) = \rho_*$ while  $\rho(t_*) = 2 \rho_*$ 
  and $\rho$ is monotonically increasing  over the interval    $t_1 < t < t_2$
  with  the derivative bound
  \begin{equation}
  |\dot \rho (t)  | \le c_3. 
  \label{eqn:rhobd}
  \end{equation}
    If the solution is oscillatory,
  then  we can   take $t_*$ such that $\dot \rho (t_*) = 0$ and continuing further, 
 find   $t_3> t_*$ such that $\rho(t_3) = \rho_*$, and  
   $\rho$ decreases  monotonically  over$[t_*, t_3]$ with the bound (\ref{eqn:rhobd})
 in place. 
  \end{prop}
       In the oscillatory case,  the constant $c_3$ can be taken as small as we like.
    In the escape case  the limit  $\nu_{\infty}: = \lim_{t \to \infty} \dot \rho (t)$
    exists and we can take for $c_3$ any constant greater than $\nu_{\infty}$.(See the  (\ref{eqn:rhoDE2}) of the appendix regarding this limit.) 
    
  We      show how  the bounds
 of the proposition imply the    desired bound   (\ref{eqn:interval}) on $\ell$.  
  We measure the length $\ell$
  of the domain of the arc of solution guaranteed by the   proposition  
\[
\begin{array}{lcl}  
 \ell  &= & \int d \sigma \\
  &= &  \int{{d \sigma} \over {dt}} dt \\
   &= & \int { dt  \over {R^2 \lambda }} \\
     & \ge  & K  \int_{t_1} ^{t_*} {dt \over {\rho ^2}}  \\
      & = & K \int_{t_1} ^{t_*} {dt \over {d \rho}} {{d \rho} \over {\rho^2}} \\
     & \ge  &{K \over c_3}  \int_{\rho_*} ^{2 \rho_*} {{d \rho} \over {\rho^2}} \\
     & = & { K \over  2c_3} { 1 \over {\rho_*}}  \\ 
      & \ge &   {C \over R_*}    
        \end{array}
     \]
  which is the desired bound.  
   
   \subsection{Proving   Proposition \ref{prop:longsolutions} }
  The proof divides into two cases,   escape and oscillatory.
  Both cases rely on   the  
     inequality:  
    \begin{equation}
       \label{eqn:rhoDE2}
-c_- / \rho^2 <   \ddot \rho < - c_+/ \rho^2, 
    \end{equation} 
   valid for   $\rho > \rho_0$ with $\rho_0$ large enough.
    As usual, the constants  $c_- > c_+ > 0$ depend only on   the masses.  
By taking $\rho_0$ arbitrarily large, we can make
    $c_-, c_+$ arbitrarily close to each other
    and to the total mass. See the appendix, lemma \ref{lem:asympKepler} for
    the proof of (\ref{eqn:rhoDE2}).      
    
{\bf Case 1: Escape.} Say that $\rho(t) \to \infty$ with $t$.
According to (\ref{eqn:rhoDE2}) its speed $\dot \rho$ 
decreases with increasing $\rho$.  For $t$ sufficiently
large, we have $\dot \rho (t) >0$, for otherwise
we would have arbitrarily large times at which $t$
 would turn back around and $\rho$ would 
decrease, contradicting escape. (See the comparison lemma
immediately below, and the appendix for more details. ) 
  It follows that $\dot \rho (t)$ is monotonically decreasing with increasing $t$
 and so  tends to  a    limit
$\nu_{\infty} \ge 0$.  Given any $\epsilon > 0$, choose $t_*$ large enough so that 
$0  < \dot \rho \le \nu_{\infty} + \epsilon$ while $\rho(t_*) := \rho_* > \rho_0$.   Then
for all $t > t_*$  we have $\rho(t) \ge \rho_*$ while  
$0 < \dot \rho (t) \le  \nu_{\infty} + \epsilon$.  
Thus $\rho$ travels between $\rho_*$ and $\infty$ all the
while satisfying $|\dot \rho| \le  c_3 = \nu_{\infty} + \epsilon$.

{\bf Case 2:  Oscillatory.} We     use inequality   (\ref{eqn:rhoDE2}) 
in conjunction with: 
      \begin{lemma} \label{lemma:cflemma}
     [Comparison Lemma]. 
 Consider   three  scalar  differential equations  
$\ddot x_- = F_- (x_-)$,   $\ddot x = F (x, t)$, $\ddot x_+ = F_+ (x_+)$
with $C^1$ right hand sides 
  satisfying   $F_- (x) < F (x, t) < F_+ (x) <0 $ for $x > x_c$, $x_c$ a fixed constant.
Suppose that    $F_-(x)$ and $F_+(x)$  are  monotone   increasing for $x > x_c$   
Let
$x_-  (t)$,  $x (t) $, $x_+ (t)$ be the solutions to   their   respective differential equation
  sharing `at rest'   initial conditions at time $0$:  
$x_- (0) = x_1 (0) = x_+ (0) : = x_*  >x_c ,  \dot x_- (0) = \dot x (0) =  \dot x_+ (0) =0$.
 %Write  $t_-(x), t(x), t_+ (x) > 0$ for the time it takes 
% the solution $x_{\pm} (t), x(t)$
% to reach the point $x$.  (If $x_i$ never reaches $x$ then $t_i (x) = + \infty$).
 Then,  for all times $t$ such that $x_- (t) \ge x_c$
 we have
 %as long as all functions exist and are positive  we have
 \begin{itemize}
\item{ (1)}$x_- (t) \le  x (t) \le x_+ (t) \le x_*$ with equality only at $t=0$,  and
%\item{ (2)}  $t_- (x) \le t(x) \le  t_+ (x)$ with equality only for $x = x_0$.  
\item{(2) } $dx_-  (t) /dt  <  dx(t)/dt < dx_+ (t) /dt $ for $t >0$ 
and  $dx_-  (t) /dt >  dx (t)/dt  >  dx_+ (t)/dt $ for $t < 0$
  \end{itemize}
   See Figure 1.   \end{lemma}

\begin{figure}
\centering
\includegraphics[scale=0.3]{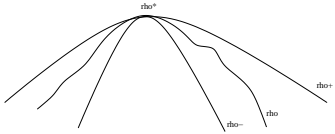}
\caption{Figure 1. Functions in the comparison lemma.}
\label{fig:rho}
\end{figure}

 We   prove the lemma below. 
Assuming the lemma we   continue with the
 proof of the proposition in the   oscillatory case.  Let $\rho_0$ large be chosen so that
 the  estimate of (\ref{eqn:rhoDE2}) is in force for $\rho > \rho_0$,  with the constants $c_+, c_-$  sufficiently close to each other.
 How close is detailed in the next paragraph.  
Since the solution is    oscillatory,  
given any $\rho_* > 0$  we can find
times $t_*$ 
arbitrarily large such that 
$$\rho (t_*) : = 2 \rho_* $$
and 
$$\dot \rho (t_*) = 0.$$
Since $\rho_*$ is arbitrarily large, we may suppose that  
 with ${1 \over 2} \rho_* \ge \rho_0$.  
 The comparison lemma sandwiches  $\rho$ between the solutions $\rho_+, \rho_-$
to the `bounding' differential equations:
$\ddot \rho_{\pm} = -c_{\pm} / \rho_{\pm} ^2$ of (\ref{eqn:rhoDE2})
which share initial conditions with $\rho$ at $t = t_*$.  
See figure 1.  Thus  
$$\rho_- (t) < \rho(t) < \rho_+ (t);  t \in I$$
where $I$ is an interval containing $t_*$
  such that  for $t \in I$ the bound $\rho (t) > \rho_0$  
  needed to obtain (\ref{eqn:rhoDE2}) is in force.

We can describe the comparison solutions $\rho_{\pm}$ in sufficient detail
by   using
the scaling symmetry of Kepler's equation.  
 Let $\phi(t)$ be the solution to the model Kepler equation
 $\ddot \phi = -1/\phi^2$ with initial conditon
 $\phi(0) = 1, \dot \phi(0) = 0$.
 Then  
 $$\rho_+ (t) = \lambda \phi ( \lambda^{-3/2}  \sqrt{c_+} (t - t_*))$$
 and
 $$\rho_- (t)= \lambda \phi ( \lambda^{-3/2}  \sqrt{c_-} (t - t_*))$$
 where we take
 $$\lambda = 2 \rho_*$$
 to guarantee agreement of initial conditions at $t = t_*$. 
By taking $\rho_0$ sufficiently large we can make $c_-$
 arbitrarily close to $c_+$. Consequently, for $\rho_0$ large enough
 we will have that $1/4 \le \phi(\sqrt{c_-} \tau_f)$
 where $\tau_f >0 $ is time such that 
 $\phi(\sqrt{c_+} \tau_f) = 1/2$.
 Now the times $\tau$ and $t$ for the scaled solutions are related
 by  $\tau = \lambda^{-3/2} (t - t_*)$. It follows at the time $t_2$ corresponding to
 $\tau_f$ we have  
 $$\rho_+ (t_f ) = \rho_* $$
 and $\rho_0 < \rho_*/2 <  \rho_- (t_f) < \rho(t_f) \le  \rho_*$.
Over the time interval $[0,  \tau_f]$   the   uniform 
 derivative bound 
 $- k < \sqrt{c_-} {{d \phi} \over {d \tau}}  (\sqrt{c_-} \tau)  < 0$ 
 holds.  Under  the  scaling and translation  symmetry used to make $\rho_{\pm}$
 we find that   velocities
 transform by 
 $v(t)  =  {1 \over {\sqrt{\lambda}}} 
 {{d \phi} \over {d \tau}} (\lambda^{-3/2} (t-t_* ))$.
 Consequently $ k/\sqrt{ \lambda} < \dot \rho_- < 0$
 during the time invterval $[t_*, t_f] $.  
 By the comparison lemma then
 $$ -k/ \sqrt{ 2 \rho_*} < \dot \rho < 0$$
 over this same time interval.  Now $\rho(t_f)$ may be less than
 $\rho_*$ but $\rho$ is monotonically decreasing.
 So   take for $t_3$ the unique time in the interval $[t_*, t_f]$
 such that   $\rho(t_3) = \rho_*$. This completes
 the  argument in the   oscillating case   for the decreasing interval $[t_*, t_3]$  of $\rho$.
 The argument for the  increasing arc $[t_2, t_*]$ of $\rho$ is   the time reversal
 (about $t_*$)  of this argument.  
 
 QED (for the proposition.)
 
\smallskip

  {\bf Proof of Lemma.}
  
 Proof of lemma.   (1) follows from (2) by integration.
% (3)    follows (2) and the negativity of $F_+$.   Since $F_+$ is negative, $dx_+/dt$ is %everywhere
% negative for $t  > 0$, and so $x_+$ is
% monotone decreasing.   Thus there is a 
% finite unique $tg_1$ such that $x_1 (t_1) = x$.
% Suppose $t_1 = t_1 (x)$ so that
% $x_1 (t_1) = x$.  And by (2) we   have for $0 \le s \le t_1$ that
 % $x  < x_1 (s) < x_2 (s)$ so it is impossible to have
 %$x_2 (s) = x$ for   $s$ in this interval. That is, $t_2 (x) > t_1 (x)$.
 We will just  prove (2) in the $-$ case, i.e.  the inequality $dx_- /dt <  dx/dt$ for $t > 0$,
 The argument in the other cases is  identical. 
Looking at the Taylor expansions  of $x_-, x$ at $t = 0$ we see that 
 the  inequality     holds in  a small right-hand neighborhood
 of $0$, say  $(0, \delta)$.  Now   proceed by contradiction.
  If the inequality fails before $x_-$ reaches $x_c$,
 then there is a   $t$  with
 $dx_- (t) /dt  \ge dx (t) /dt $.   Let $t_*$ be the first such $t > 0$
 such that $dx_-  (t) /dt  = dx /dt$.  We have $t_* > \delta$.
By integration, $dx_- (t_*)/dt  = \int_0 ^{t_*} F_- (x_1 (t)) dt$
 and $dx (t_*)/dt  = \int_0 ^{t_*} F (x(t) , t)dt$.  These two 
 integrals are equal.  But $F (x(t), t)) > F_- (x(t))$.
 And in the interval $(0, t_*)$ we have $dx_- /dt < dx /dt$, and
 so, by integration $x_- (t) < x (t)$.  Then  $F_- (x (t)) > F_- (x_- (t))$
  by   the monotonicity of $F_-$.
 So $ \int_0 ^{t_*} F (x (t) , t)dt >  \int_0 ^{t_*} F_- (x_- (t)) dt$, 
 contradicting the   equality of the two integrals.  
 
 QED
 
 \medskip 
 {\bf Remark.}  Differential  inequalities
involving  $\rho$ are much better behaved at large $R$ (and hence $\rho$)
than those involving  $R$.  The $R$ differential equation is   the Lagrange-Jacobi identity 
   $2 \ddt {  R \dot R} = 4 H + 2U$ and yields a huge second derivative for $R$
   when $r = r_{12}$ is sufficiently small.    Thus $R$ can
      oscillates wildly, despite the fact
   that the bound (\ref{eqn:Rbd}) is in force.

\medskip

{\bf Discussion. Open questions.} 
  Could the theorem   hold for arbitrary energy $H$
and angular momentum $J$? No.   It   does not hold for $H > 0$ and  $J = 0$.
For the    direct method of
the calculus of variations    yields action-minimizing  hyperbolic escape orbits
which leave  triple collision and tend to any desired noncollinear point
of the shape sphere in infinite time.  The reflection argument
(see, eg.  \cite{MontChenc}) shows that these minimizers
never become collinear.    The  theorem might hold
for $H = J = 0 $  but I suspect not.  
In this case there is a manifold of    parabolic escape orbits whose
shapes tend to  Lagrange.  I would guess some of these have no 
 syzygies, but this is just a guess.   
For $H < 0$,  and general $J \ne 0$
the theorem is   false,  at least for  mass ratios in which one mass dominates.  
 For in this case,  the  near-circular 
  Lagrange solutions are KAM
stable, and so are surrounded by a  nearby  cloud
of KAM torii on which the solutions  stay near Lagrange,
and hence away from $z = 0$ for all time.    It is possible
that for some values of $H<0$, $J \ne 0$ and some values of 
the   mass ratios that the theorem continues to hold.
If the   Dziobek constants $J^2 H$ and mass ratios are such that
   the Lagrange solution is 
  unstable (which is the case for  nearly equal masses and 
  $J^2 H$ being a value near that which supports the circular
  Lagrange solution) then there is some chance for  the theorem
   to  hold.   
  
  According to the theorem, 
  all   solutions bi-asymptotic to triple collision   except
  for the  Lagrange solution  have syzygies.
  This number is necessarily finite.   What numbers are possible? 
    Is any finite number of syzygies achieved?   Is any finite syzygy sequence
  realized?   Write $m(x)$ for the time interval on which the solution $x$ is defined.
(Thus $m(x) = + \infty$ for all solutions except those bi-asymptotic to
  triple collision.) Is it true that $m(x)$ is minimized (among all solutions $x$  with 
  $J = 0$ and $H < 0$ fixed)  by the Lagrange solution?

\section{Appendix: Bounds  near $\infty$ for negative energy.}

We suppose the total energy $H = K/2 - U$
to be negative and write $H = -|H|$.   Then
\begin{equation}
\label{eqn:Ubd1}
 U \ge |H|
\end{equation}
Write $r = min \{r_{ij}: i \ne j \}$
for the minimum of the intermass distances.
Then there is a   constant $c$ depending only on the masses
such that
\begin{equation}
\label{eqn:Ubd2}
c/r \ge U.
\end{equation}
(For instance, if the masses are all equal to $m$ then $c = 3m^2$.)
It follows that 
$$c/r \ge  |H|$$
or
\begin{equation}
\label{eqn:rbd}
c/|H|  \ge r.
\end{equation}

Let us suppose that $12$ realize the minimum distance:
$$r = r_{12}.$$
Associated to the decomposition $12; 3$
we have  Jacobi vectors and their lengths:  
$$\zeta = x_1 - x_2  \hskip .5cm,   |\zeta| =r  $$
$$\xi = x_3 - x^{12}_{cm}   \hskip .5cm,   |\xi| =\rho.$$
Here 
\[
\begin{array}{lcl}  
 x^{12}_{cm}  & :=  &(m_1 x_1 + m_2 x_2 )/(m_1 + m_2 )  \\ 
 & : = & \mu_1 x_1 + \mu_2 x_2
\end{array}
 \]
is the $12$ center of mass and 
$$\mu_1 = m_1/ (m_1 + m_2)  \hskip .5cm,  \mu_2 = m_2 /(m_1 + m_2).$$
 One computes  
\begin{equation}
\label{eqn:R2} 
  R^2 = \alpha_1 r^2 +  \alpha_3 \rho^2
    \end{equation}
  where 
 \begin{equation}
\label{eqn:redmasses}
  \alpha_1 =m_1 m_2 / (m_1 + m_2), \alpha_3 = (m_1 + m_2) m_3/(m_1 + m_2 + m_3).
    \end{equation}
 from which it follows that
   \begin{equation}
\label{eqn:rhoRbds}
   c_a \rho \le R \le C_a \rho 
   \end{equation}
  and 
  \begin{equation}
\label{eqn:Ubd3}
   U \ge C/ \rho
  \end{equation}
 
 Set $\hat U = RU$. Combine (\ref{eqn:Ubd1}), (\ref{eqn:Ubd2}), (\ref{eqn:Ubd3})
  with (\ref{eqn:rhoRbds}) to get that
  $C \rho /r \ge \hat U \ge R|H|$ or
  $$ {C \over  R|H|} \ge r/\rho$$
  which asserts that by making $R$ or $\rho$  large we can make the ratio
  $r/\rho$ as small  as we wish. We view  $r/\rho$ as a  perturbation parameter.
From the last inequality   it follows that  
  for every $\epsilon > 0$ there is a $\rho_0$  (or  $R_0$) sufficiently large 
  such that
  $\rho \ge \rho_0$ ($R \ge R_0$) implies that $r/\rho < \epsilon$.
  {\bf In what follows, let $\epsilon$ small be   given, and
  suppose we have chosen the  corresponding  $\rho_0$ (or  $R_0$) be taken
  so as to guarantee $r / \rho < \epsilon$.   And let $c$, $C$, .. denote constants
  depending only on this $\rho_0$, the masses, the total energy,
  and, in a moment, the total angular momentum.  }

  We can  express the other distances in terms
  of $\xi, \zeta$:   
   \begin{equation}
\label{eqn:distances}
r_{13} = \| \xi - \mu_1 \zeta \| \hskip .13cm ; \hskip .13cm r_{23} = \| \xi + \mu_2 \zeta \| 
\end{equation}
 where $\mu_i$ are the reduced masses described above. 
  Note that $r_{13} = \rho + O(\epsilon)$ and $r_{23} = \rho + O(\epsilon)$.
 
 We have $$H  =  H_{12}  + H_3  + g \hskip .4cm (5a) $$
where
$$H_{12} = {1 \over 2} \alpha_1 \| \dot \zeta \|^2 - \beta_1/r, \hskip .4cm (5b)$$
$$H_3 = {1 \over 2} \alpha_3 \| \dot \xi \|^2 -  \beta_3/\rho, \hskip .4cm (5c)$$
where $\alpha_1, \alpha_2$ are given in eq. (\ref{eqn:redmasses} )
 and where  
$$\beta_1 = m_1 m_2,   \beta_3 = (m_1 + m_2) m_3$$ 
and where the ``error term''  $g$ is given  by
\[
\begin{array}{lcl} g & = & {{(m_1 + m_2) m_3} \over{\|\xi \|}} - {{m_1 m_3} \over {\| \xi - \mu_1 \zeta \|}} - 
{{m_2 m_3} \over { \| \xi + \mu_2 \zeta \|}} \\
 & =  &(  m_1 m_3 \mu_1  - m_2 m_3 \mu_2 ) {1 \over { \|\xi \|^3}} \langle \xi, \zeta \rangle  
 + O({1 \over {\rho^3}}) \\
 \end{array}
\]
Note that
\begin{equation}
\label{eqn:gbd1}
 |g| \le C \epsilon /\rho.
\end{equation}
We will also need bounds for the gradients of $g$:
$$ g_{\xi} =   c \zeta / \rho^3 + O (\rho^{-4})$$
so that
\begin{equation}
\label{eqn:gbd2}
| g_{\xi} |\le C \epsilon/ \rho^2
\end{equation}

 If we set $g = 0$ then $H$ becomes the Hamiltonian for  two uncoupled
Kepler systems.    The   next lemma describes some details of
the asymptotics of this decoupling  as $\rho \to \infty$.      Introduce 
 $$J_{12} = \mu_1 \zeta \wedge \dot \zeta,$$
 the angular momentum  of the $12$ system,
  and the radial and transverse separation velocities
 $\nu, V^{\perp}$ by  
 $$\dot \xi = \nu \hat \xi + V ^{\perp} \hskip .3cm ;   \hskip .3cm \dot \rho = \nu$$
 where $\hat \xi = \xi /\rho$ is the unit vector in the $\xi$
 direction and where $V^{\perp}$ is orthogonal to $\xi$. 
  
 \begin{lemma} 
 \label{lem:asympKepler}
   Consider any solution to the three-body problem
   along which $\rho (t) \ge \rho_0$ with $\rho_0$ as above. 
There exists a positive constant $c$,
 depending only on  the total energy $H$
  and  total angular momentum $J$,  the masses  $m_i$
   and $\rho_0$ such that
  \begin{itemize} 
 \item{(a)} $|J_{12}| \le c$
 \item{ (b)}  $\|V ^{\perp} \| \le c/\rho$
  \item{(c)} $|\ddot \rho + M/ \rho^2|   \le \epsilon/ \rho^3$ where
 $M =  m_1 + m_2 + m_3$ is the total mass. 
  \end{itemize}
 \end{lemma}

  Similar bounds hold
  for the time derivatives of  $H_{12}$,
  $H_3$,  $J_3$,  $\hat \xi$, and the $12$ and $3$   Laplace (or Runge-Lenz) vectors.
   
\smallskip
{\bf  Proof of  Estimate (a):}  We   show that 
  $\|J_{12} \| \le  \beta_1 ^2/|H| + O(1/\rho_0)$.
  We have   $J_{12} = \alpha_1 \zeta \wedge \dot \zeta$, $ \| \zeta \wedge \dot \zeta \|^2 \le  \|\zeta\|^2 \|\dot \zeta \|^2$ and   $|\zeta |^2 = r^2$.
  It follows that  
  $$|J_{12}| \le \alpha_1 r^2   \|\dot \zeta \|^2 . $$   
  Set $H^{\prime} =  H_3 + g $.
  Note that $-H^{\prime} \le \beta_3/ \rho + \epsilon/\rho \le c/\rho_0$,
  and that 
 $ H_{12} = H -  H^{\prime}$.
 It follows that 
 But $H_{12}  \le  -|H| + c/\rho_0$.
 Use the formula for   $H_{12}$ to rewrite this
 inequality as    
 $$\alpha_1 \| \dot \zeta \|^2 \le -2|H| + 2 c/\rho_0 + 2 \beta_1/r$$
 Multiply through by $r^2$  to get
 $$|J_{12}| \le \alpha_1 ^2 r^2 \| \dot \zeta \|^2 \le 
   [ -2|H| + 2 c/\rho_0 ] r^2 +  2 \beta_1  r$$
The right hand side  is a quadratic function of $r$ with negative quadratic term.
The   maximum value of this quadratic function  is   $({1 \over 2})(2 \beta_1)^2/(2|H| - 2c/\rho_0))
= \beta_1 ^2/|H| + O(1/\rho_0)$.  Thus
$|J_{12}| \le \beta_1 ^2/|H| + O(1/\rho_0)$.

\smallskip
{\bf Proof of estimate (b). } 
We have $\alpha_3 \xi \wedge \dot \xi = \alpha_3 \xi \wedge V ^{\perp}$,
$|\xi \wedge V^{\perp}| = \rho |V^{\perp}|$ 
and $\alpha_3 \xi \wedge \dot \xi  = J - J_{12}$.
Thus 
 \[
\begin{array}{lcl}
 \alpha_3 \rho |V^{\perp}| &  =  &| J - J_{12}|
\\
& \le & |J| + |J_{12}|
  \end{array}
\]
Now use estimate (a).

\smallskip
{\bf Proof of  estimate (c).} 
We have
 $\dot \rho : = \langle \hat \xi, \dot \xi \rangle.$
so that
\begin{equation}
\label{eqn:c1}
\ddot \rho =  \langle  {d \over dt}  \hat \xi, \dot \xi \rangle +  \langle \hat \xi, \ddot \xi \rangle.
\end{equation}
We compute that 
$\langle {d \over dt} \hat \xi, \dot \xi \rangle = -  {{\dot \rho^2} \over \rho} + {{\| \dot \xi^2 \|} \over \rho} = \|V^{\perp} \|^2/ \rho$
so that by estimate (b): 
\begin{equation}
\label{eqn:c2}
| \langle  {d \over dt}  \hat \xi, \dot \xi \rangle | \le c/ \rho^3
\end{equation}
Now use Newton's equations for $\xi$ 
$$\alpha_3 \ddot \xi =   U_{\xi} = - \beta_3 \xi/ \rho^3 + g_{\xi}.$$ 
which yields 
$$\ddot \xi  = - M \xi /\rho^3 +{1 \over \alpha_3} g_{\xi}$$
because  $\beta_3/\alpha_3 = M$.
Thus 
\begin{equation}
\label{eqn:c3}
  \langle \hat \xi, \ddot \xi \rangle = - M/ \rho ^2 + {1 \over \alpha_3} \langle \hat \xi, g_{\xi} \rangle.
  \end{equation}
Using  the estimate (\ref{eqn:gbd2}) on the gradient
$g_{\xi}$ of $g$, with 
equations (\ref{eqn:c1}), \ref{eqn:c2}), and (\ref{eqn:c3})
we get the 
desired result
$| \ddot \rho  + M/ \rho^2 | \le c/ \rho^3$.

QED

As a general   reference for some  of the inequalities
appearing  here, and many others, see  Marchal, \cite{Marchal}
esp. pp. 327-7 , equations (885)-(894).

\end{document}